\documentclass[11pt]{article}
\usepackage{amsfonts}

\usepackage{amssymb}
\usepackage{latexsym}
\usepackage{amscd}
\usepackage[centertags]{amsmath}
\usepackage{xypic}

\newtheorem{thm}{Theorem}[section]

\newtheorem{lem}[thm]{Lemma}
\newtheorem{rem}{Remark}

\numberwithin{equation}{section}

\newcommand{\calF}{\mathcal{F}}
\newcommand{\calG}{\mathcal{G}}

\newcommand{\calS}{\mathcal{S}}

\newcommand{\Enum}{\mathbb{E}}
\newcommand{\Pnum}{\mathbb{P}}

\newcommand{\Znum}{\mathbb{Z}}
\newcommand{\Nnum}{\mathbb{N}}

\newcommand{\Om}{\Omega}
\newcommand{\om}{\omega}

\newcommand{\qed}{\hfill $\Box$}

\newcommand{\supp}{\rm supp}

\newcommand{\Proof}{ \noindent {\bf Proof.}\ }

\newcommand{\beqn}{\begin{eqnarray}}             
\newcommand{\eeqn}{\end{eqnarray}}               
\newcommand{\beq}{\begin{eqnarray*}}             
\newcommand{\eeq}{\end{eqnarray*}}               
\begin{document}
\title{Average Entropy of the Ranges for Simple Random Walks on Discrete Groups}
\author{Xin-Xing Chen$^{1,4}$ \quad Jian-Sheng Xie$^{2,4}$ \quad Min-Zhi Zhao$^{3,4}$\\
%
{\footnotesize 1. School of Mathematical Sciences, Shanghai Jiaotong University, Shanghai 200240,  China,}\\
{\footnotesize E-mail: chenxinx@sjtu.edu.cn}\\
 {\footnotesize 2. School of Mathematical Sciences, Fudan
University, Shanghai 200433,   China,}\\
 {\footnotesize E-mail:jsxie@fudan.edu.cn }\\
{\footnotesize 3. School of Mathematical Sciences, Zhejiang
University, Hangzhou 310027,   China,}\\
 {\footnotesize E-mail:zhaomz@zju.edu.cn }\\
 {\footnotesize 4. Shanghai Center of Mathematics, Fudan University,
Shanghai 200433,   China.}} 
\date{}
\maketitle
\begin{abstract}
Inspired by Benjamini et al (Ann. Inst. H. Poincar\'{e}  Probab. Stat. 2010) and Windisch
(Electron. J. Probab. 2010), we consider the entropy of the random walk range formed by a
simple random walk on a discrete group. It is shown in this setting the existence of a quantity
which we call the average entropy of the ranges. Some equivalent conditions for the vanishing of
the average entropy of the ranges are given. Particularly, the average entropy of the ranges vanishes
if and only if the random walk is recurrent or escaping to  negative infinity  without left jump.
In order to characterize the recurrence further, we study the average entropy of the weighted
digraphs formed by the random walk. We show that the random walk is recurrent if and only if the
average entropy of the weighted digraphs vanishes.

\end{abstract}

{  \verb"Keywords":}  random walk,  entropy, range, recurrent

{  \verb"AMS 2000 subject classifications":}   60G50,60J10,05C25


\section{Introduction}\label{sec: 0}
Let $R_n$ be the random walk range formed by the first $n$ steps of a  random walk. The study
of $R_n$ seems to be initiated by Dvoretzky and Erd\"{o}s   \cite{D&E1950}, where they
established strong law of large numbers for the size of the random walk range (i.e. $|R_n|$)
in lattice cases $\Znum^d$ with $d \geq 2$. Probabilistic scientists also found that such
quantity satisfies the central limit theorems \cite{J&P1971} \cite{J&P1972'} \cite{J&P1974}
\cite{LeGall1986} and laws of iterated logarithm \cite{B&K2002} \cite{J&P1972}. Many other
quantities about the structure of $R_n$, e.g. the radius of the largest covered ball, were
discussed in Erd\"{o}s and Taylor \cite{E&T1960/A,E&T1960/B} and R\'{e}v\'{e}sz \cite{Revesz2005}.
Recently, Benjamini et al \cite{BKYY2010} and Windisch \cite{Windisch2010} discussed the order
of $H (R_n)$, the entropy of random walk range for lattice cases and graph cases (with the
symmetrical distributions) respectively. It is natural to ask what happens for the general case
of simple random walks on discrete groups.

Inspired by their works, and bearing in mind the question above, we continue to research on
the entropy of random walk range for a simple random walk on a discrete group. In this setting
we find that the limit $\lim\limits_{n\to\infty} \frac{H (R_n)}{n}$ always exists (which may
be known to the authors of \cite{BKYY2010, Windisch2010}) and we will call it the average
entropy of the ranges. Furthermore, we give different equivalent conditions for the average
entropy attaining its critical values (the possible minimum value 0 and the possible maximum
value). In the meanwhile we also relate the tail $\sigma$ algebra with the vanishing of the
average entropy of the ranges following a work of Ka\u{\i}manovich and Vershik \cite{Kaimanovich1983}.
An important consequence is that, when $0<H(R_1)<\infty$, the random walk is recurrent if and
only if the tail $\sigma$-algebra of $R_n$ is trivial. As is seen, the vanishing of the so-called average
entropy of the random walk ranges cannot be used to characterize recurrent property which leads
us to consider the average entropy of the weighted digraphs for random walks.

\section{Settings and Main Results}\label{sec: 1}
Let $G$ be a countable group with the identity element $e$. Let $\mu$ be a probability measure
on $G$ such that $G$ is generated by the support of $\mu$. Let $(\Omega, \mathbb{P})$ be the
respective infinite product space, i.e.,  $\Om=G^{\Nnum}$ and $\Pnum=\mu^{\Nnum}$. Let $\theta$
be the natural time shift map on $\Omega$, i.e.,
$$
\theta ((\omega_1,\omega_2,\cdots))=(\omega_2,\omega_3,\cdots).
$$
Then the system $(\Omega,\mathbb P;\theta)$ is measure-preserving and ergodic. Next, we will
define a random walk. For $n\ge 1$ and $\om=(\om_1,\om_2,\cdots) \in \Omega$, set $X_n (\om)
=\om_n$. Then $X_1,X_2,\cdots$ are independent and have the common  law $\mu$. Let $S_0=e$, and
set
$$
S_n=S_{n-1} X_n~~{\rm for~}n\ge 1.
$$
Then $S=(S_0,S_1,S_2,\cdots)$ is a simple random walk on $G$ which is called a $(G, \mu)-$random
walk throughout this paper. 

Let $R_n=\{g\in G:\exists 0\le i\le n, S_i=g\}$, the range of $S$ until time $n$. Let $\Gamma_n
=(R_n, \overrightarrow{E}_n, W_n)$ be the random weighted digraph generated by the first $n$ step,
where
$$
\overrightarrow{E}_n=\{(x,y): \exists i< n, S_{i}=x, S_{i+1}=y\}
$$
and
$$
W_n(\overrightarrow{a})=\sum_{i=0}^{n-1} 1_{\{ S_i=x,
S_{i+1}=y\}}{\rm~~~for~}\overrightarrow{a}=(x,y) \in \overrightarrow{E}_n.
$$
So the weight function $W_n(\overrightarrow{a})$  records the number of jumps passing through the
edge $\overrightarrow{a}$ before time $n$.

As usual, we write $H(Y)$ for the entropy of a discrete random element $Y$, that is
$$
H(Y)=-\sum_y \Pnum(Y=y)\ln \Pnum(Y=y)
$$
with usual convention $0\ln 0=0$; it is called Shannon entropy with reference to Shannon's work
\cite{Shannon1948}. In this paper we discuss mainly on the limits of $\frac{H(R_n)}{n}$
and $\frac{H(\Gamma_n)}{n}$. To avoid triviality, we assume $0<H(X_1)<\infty$ further. For each
$A\subset G$ and $x\in G$, we set
\begin{equation}\label{eq: 1.1}
p_n(A)=\Pnum (R_n=A)~~~~~{\rm and~~~~~} q_n(A,x)=\Pnum (R_n=A, S_n=x).
\end{equation}
Then we have our first result.
\begin{thm}\label{thm: 1.1}
There is a constant $0\le h\le H(X_1)$ such that
\begin{equation}\label{exsitence1}
\lim_{n\to\infty}\frac {H(R_n)}{n}=\lim_{n\to\infty} \frac{H(R_n,S_n)}{n}=h
\end{equation}
and
\begin{equation}\label{exsitence2}
\lim_{n\to\infty}\frac {-\ln p_n(R_n)}{n}=\lim_{n\to\infty} \frac{-\ln q_n(R_n,S_n)}{n}=h\,\,\, a.s.
\end{equation}
The above statements (\ref{exsitence1}) and (\ref{exsitence2}) also hold when $R_n$ is replaced by $\Gamma_n$,
while the constant $h$ may be different.
\end{thm}

Write  $h_R(\mu):=\lim \limits_{n\to\infty}\frac {H(R_n)}{n}$ and call it the average entropy
of the range of a $(G, \mu)$-random walk. Similarly,   write $h_{\Gamma} (\mu):=\lim \limits_{n\to\infty}
\tfrac{H(\Gamma_n)}{n}$ and call it the average entropy of the traces. See, e.g., \cite{RS1982} \cite{BGL2007}
\cite{BM2012} \cite{Shiraishi2010} \cite{Shiraishi2012} for related studies on random walk traces.

The average entropies have close relations with tail $\sigma$-algebras. Write $R_\infty=\cup_n R_n$, and put
\begin{equation}\label{e:1120w}
\calG_n :=\sigma (R_k, S_k: k \geq n),~~ \calG_\infty:= \bigcap_{n \geq 0} \calG_n.
\end{equation}
In other words,  $\calG_\infty$  is the tail $\sigma$-algebra of $(R_n, S_n)$. A well known result (see
\cite{Kaimanovich1983} or \cite[Chap. 12]{Benjamini}) indicates the existence of the limit $h_S (\mu)
:=\lim \limits_{n \to \infty} \frac{1}{n} H (S_n)$ and
\begin{equation}\label{eq:S-entropy-formula}
h_S (\mu)=\frac{1}{k}[ H (X_1, \cdots, X_k) -H (X_1, \cdots, X_k|
\calS_\infty)], \quad \forall k \geq 1,
\end{equation}
where $\calS_\infty$ is the tail $\sigma$-algebra of $S_n$. In a similar spirit, we have the following
formula for $h_R(\mu)$.
\begin{thm}\label{thm: 1.4'}
For any $k \geq 1$,
\begin{equation}\label{e:1250}
h_R (\mu) =\frac{1}{k}[H (X_1, \cdots, X_k) -H (X_1, \cdots, X_k, R_\infty
\circ \theta^k ~|~ \calG_\infty)].
\end{equation}
\end{thm}

It is well known that (see \cite[Chap. 12]{Benjamini}) $h_S (\mu)=0$ if and only if $\calS_\infty$ is trivial.
Inspired by this fact, we would go further to find out equivalent conditions for the vanishing of the average
entropy of the range $h_R(\mu)$. As we know, the symmetric simple random walk on $\Znum^d$ with $d\le 2$ is
recurrent and $R_\infty=\Znum^d$ almost surely. A result in \cite{BKYY2010} tells us $H (R_n)=O (\log n)$ for
$d=1$ and $H (R_n)=O (\frac{n}{\log n})$ for $d=2$; Thus the average entropies of the ranges in these two cases
are zero. It turns out that recurrent random walks always make $h_R(\mu)$ vanish; but still there is another
special type of transient random walks bearing this property. To state this result in the simplest form, we
need more notations. Set
$$
L_\mu=\{a\in G: {\supp}(\mu)=\{a^{-1},e,a,a^2,\cdots\} {\rm~and~}
\sum_{i\ge  -1} i\mu(\{a^{i}\})<0\}.
$$
If $L_\mu\not=\emptyset$, then $G$ is isomorphic to $\mathbb{Z}$. So we may call such $S$ {\it a random walk
escaping to $-\infty$ without left jump} if $L_\mu\not=\emptyset$; clearly such a random walk is transient.


\begin{thm}\label{thm: 1.2}
The following four statements are equivalent:
\begin{enumerate}
\item[{\rm (1)}] $h_R(\mu)=0$;

\item[{\rm (2)}] $\Pnum(\exists n\ge 1, ~S_n=x)=1$ for some $x\in G$;

\item[{\rm (3)}]   $S$ is either recurrent or a random walk escaping to  $-\infty$ without left jump;

\item[{\rm (4)}] $H({R}_\infty\circ \theta|{R}_\infty)= H({R}_\infty|{R}_\infty\circ\theta)$ and
$\mathcal{G}_\infty\equiv\sigma({R}_\infty) \mod  \Pnum.$
\end{enumerate}
\end{thm}

We see that the vanishing of $h_R(\mu)$  can not  characterize the recurrence. Luckily, the following theorem
tells us that the vanishing of $h_{\Gamma} (\mu):=\lim \limits_{n\to\infty}\tfrac{H(\Gamma_n)}{n}$ can
serve this role.

\begin{thm}\label{thm: 1.3}
Process $S$ is recurrent  if and only if $h_{\Gamma} (\mu)=0$.
\end{thm}

Similar as the proof of case 1 of  Theorem 6.7.3 in \cite{Durrett2010}, we can show that the tail
$\sigma$-algebra of $(\Gamma_n,S_n)$ is trivial when $S$ is recurrent. In view of the formula (\ref{e:1250}) and Theorem \ref{thm: 1.2}, it is easy to see that  the random walk is recurrent if and only if the $\sigma$-algebra
 $\calG_\infty$ is trivial, which is also equivalent to the
triviality of the tail $\sigma$-algebra of $(\Gamma_n,S_n)$.
This  justifies
the significance of the current research.

%

The rest of the paper is organized as follows. The proofs of Theorems \ref{thm: 1.1}--\ref{thm: 1.3}
are presented in Sections \ref{sec: 2}--\ref{sec: 5} separately. In section \ref{sec: 6}, some related results are presented.
We show that $h_R (\mu)=H (X_1)$ if and only if
the escape rate of  $S$  is $1$.
We also  show that $H (R_\infty)<\infty \Leftrightarrow \Enum [|X_1| \ln |X_1|]<\infty$ for any $(\mathbb{Z},\mu)-$random walk escaping to $-\infty$ without left jump.

\section{Proof of Theorem \ref{thm: 1.1}}\label{sec: 2}
Before the proof of Theorem \ref{thm: 1.1}, we show the following
lemma.
\begin{lem}\label{l:1432}
Suppose  $p_1+p_2+\cdots p_n=1$ with $p_i>0$ for all $i$. Then for any $\alpha\ge 1$,
$$
\sum_{i=1}^np_i\ln^\alpha (1/p_i)\le
(\alpha \vee \ln n)^\alpha+(\alpha-1)^\alpha.
$$
\end{lem}
\Proof
Let $f(x)=x\ln^\alpha (1/x)$.  One can easily check that $f''\le 0$ on $(0, e^{1-\alpha}]$. So
\begin{align*}
\sum_{i=1}^n p_i\ln^\alpha (1/p_i) \le & (\alpha-1)^\alpha\sum_{p_i\ge e^{1-\alpha}}p_i+ \sum_{p_i<e^{1-\alpha}}
p_i\ln^\alpha (1/p_i)\\
\le& (\alpha-1)^\alpha + \left(\sum_{p_i<e^{1-\alpha}} p_i\right) \ln^\alpha
\left(\frac{|\{i\le n: p_i<e^{1-\alpha}\}|}{\sum\limits_{p_i<e^{1-\alpha}}p_i}\right)\\
\le&(\alpha-1)^\alpha +\max_{x\in (0,1]}x \ln^\alpha( \tfrac{n\vee e^\alpha}{x}).
\end{align*}
Since $x \ln^\alpha(\tfrac{n\vee e^\alpha}{x})$ is increasing on (0,1], we finish the proof of the lemma.
\qed

\noindent\textbf{Proof of Theorem \ref{thm: 1.1}.}
Note that for any $A,B\subset G$ and $x,y\in G$, if $(R_n,S_n)=(A,x)$ and $(R_m, S_m)\circ \theta^n=(B,y)$
then
$$
(R_{n+m},S_{n+m})=(A\cup (x\cdot B), ~x\cdot y).
$$
By the Markov property,
\begin{align*}
\Pnum (R_{n+m} =A\cup(x\cdot B),~S_{n+m}=x\cdot y) \ge& \Pnum (R_n=A,S_n=x,R_m \circ \theta^n=B,
S_m\circ \theta^n=y)\\
\ge& \Pnum (R_n=A,S_n=x) \Pnum (R_m\circ \theta^n=B, S_m\circ \theta^n=y).
\end{align*}
Recall the definition of $q_n$ in (\ref{eq: 1.1}). Then the above inequality can be rewritten as
\begin{equation}\label{subadditive}
-\ln q_{n+m}(R_{n+m},S_{n+m})\le -\ln q_n(R_{n},S_{n})-\ln q_m(R_{m},S_{m})\circ \theta^n.
\end{equation}
Since $(R_{m},S_{m})\circ \theta^n$ has the same distribution with $(R_{m},S_{m})$, we take expectation on
both sides of (\ref{subadditive}) and get
$$
H(R_{n+m},S_{n+m})\le H(R_{n},S_{n})+H(R_{m},S_{m}).
$$
Furthermore, since
$$
H(R_n,S_n)\le H(X_1,X_2,\cdots, X_n)\le n H(X_1)<\infty,
$$
we can apply the Subadditive Ergodic Theorem of Kingman (see, e.g., \cite{Kingman1971}
\cite{Derriennic1983} \cite{Liggett1986}) yielding that there exists $0\le h\le H(X_1)$ such that
\begin{equation}\label{e:1649}
\lim_{n\to\infty} \frac {H(R_n,S_n)}{n}=\inf_n\frac {H(R_n,S_n)}{n}= h
\end{equation}
and
\begin{equation}\label{e:1451}
\lim_{n\to\infty} \frac {-\ln q_n(R_{n},S_{n})}{n}=h \,\,\,a.s.
\end{equation}
Moreover, since $H(R_n)\le H(R_n,S_n)=H(R_n)+H(S_n|R_n)\le H(R_n)+\ln(n+1)$, we have
$$
\lim_{n\to\infty} \frac {H(R_n)}{n}=h.
$$

Next we shall prove $\lim \limits_{n\to\infty}-\frac{\ln p_n(R_n)}{n}=h~~~ a.s..$ For $x\in A\subset G$, set
$$
p_n(x|A):=\Pnum (S_n=x|R_n=A).
$$
Then $q_n(R_n,S_n)= p_n(R_n)p_n(S_n|R_n)$. By (\ref{e:1451}), it is suffice to prove
\begin{equation}\label{e:1456}
\lim_{n \to \infty} \frac {-\ln p_n(S_n|R_n)}{n}=0 ~~\,\,a.s..
\end{equation}
Let $A\subset G$ with $\Pnum(R_n=A)>0$. Then $|A|\le n+1$. In view of Lemma \ref{l:1432},
\begin{eqnarray*}
f(A) &:=& \sum_{x\in A} \Pnum (S_n=x|R_n=A)\ln^2  \Pnum (S_n=x|R_n=A)\\
&\le& (2\vee\ln |A|)^2+1\le \ln^2 (n+1)+5.
\end{eqnarray*}
It follows immediately
\beq
\Enum[\ln^2 p_n(S_n|R_n)]&=\Enum[f(R_n)]\le \ln^2(n+1)+5.
\eeq
Hence for any $\varepsilon>0$, we have
$$
\sum_{n} \Pnum \left( \Bigl| \frac{\ln p_n(S_n|R_n)}{n} \Bigr|>\varepsilon\right)\le  \sum_{n}\frac{\ln^2 (n+1)+5}{n^2\varepsilon^2}<\infty.
$$
This proves (\ref{e:1456}) by Borel-Cantelli Lemma. Hence (\ref{exsitence1}) and (\ref{exsitence2}) hold true for $R_n$.

The result for $\Gamma_n$ can be obtained in a similar way and the details are omitted here.
\qed

\section{Proof of Theorem \ref{thm: 1.4'}}\label{sec: 3} 
\Proof
By Theorem \ref{thm: 1.1}, we have
$h_{R}(\mu)=\lim_n\frac{H(R_n,S_n)}{n}$. We also have
\begin{equation}\label{e:1243}
H (R_n, S_n) -H (R_{n-k}, S_{n-k})=H (X_1,\cdots,X_k) -H (X_1,\cdots,X_k, R_{n-k} \circ \theta^k
\big| R_n, S_n).
\end{equation}
This is because
\begin{align*}
&H (  R_{n-k} \circ \theta^k,  S_{n-k} \circ \theta^k,X_1,\cdots,X_k) \\
=& H (R_n, S_n) + H (R_{n-k} \circ \theta^k,  S_{n-k} \circ \theta^k,X_1,\cdots,X_k| R_n,S_n)\\
=& H (R_n, S_n) + H ( X_1,\cdots,X_k, R_{n-k} \circ \theta^k \big| R_n, S_n)
\end{align*}
and
\begin{align*}
H (  R_{n-k} \circ \theta^k,  S_{n-k} \circ \theta^k,X_1,\cdots,X_k)
=&H (  R_{n-k} \circ \theta^k,  S_{n-k} \circ \theta^k) +H ( X_1,\cdots,X_k) \\
 =& H (R_{n-k}, S_{n-k})+H ( X_1,\cdots,X_k).
\end{align*}
So,  if  $\lim_n H ( X_1,\cdots,X_k, R_{n-k} \circ \theta^k \big|
R_n, S_n)$ exists, then
\begin{equation}\label{e:1028}
kh_R(\mu)=H(X_1, \cdots,X_k)-\lim_n H ( X_1,\cdots,X_k, R_{n-k}
\circ \theta^k \big| R_n, S_n).
\end{equation}
Since $R_n=R_k \cup S_k (R_{n-k}\circ\theta^k)$, we have
$$
R_{n-k}\circ\theta^k=(S_k^{-1}R_n-S_k^{-1}R_k)\cup\big((R_{n-k}\circ\theta^k)\cap (S_k^{-1}R_k)\big).
$$
Thus
\begin{align*}
&H(X_1,\cdots,X_k, R_{n-k} \circ \theta^k \big| R_n, S_n)\\
=& H(X_1,\cdots,X_k, R_{n-k} \circ \theta^k \big| R_n, S_n, X_{n+j}, j=1,2,\cdots)\\
=&H(X_1,\cdots,X_k, R_{n-k} \circ \theta^k \big| \calG_n)\\
=& H(X_1,\cdots,X_k \big| \calG_n) +H(R_{n-k} \circ \theta^k \big| X_1,\cdots,X_k, \calG_n)\\
=& H(X_1,\cdots,X_k \big| \calG_n) +H( (S_k^{-1} R_k) \cap R_{n-k} \circ \theta^k
\big| X_1,\cdots,X_k, \calG_n)\\
=& H(X_1,\cdots,X_k \big| \calG_n) +H( (S_k^{-1} R_k) \cap R_{n-k} \circ \theta^k
\big| \calF_k \vee \calG_n),
\end{align*}
where $\calF_k :=\sigma (X_1, \cdots, X_k)$ stands for the $\sigma$-algebra generated by $X_1, \cdots,
X_k$. Since
$\calG_n \downarrow \calG_\infty$ we have
$$
H(X_1,\cdots,X_k \big| \calG_n) \to H(X_1,\cdots,X_k \big| \calG_\infty).
$$
So, we are left to prove  the following limit
\begin{equation}\label{e:1130}
H( (S_k^{-1} R_k) \cap (R_{n-k} \circ \theta^k) \big| \calF_k \vee \calG_n) \to H( (S_k^{-1} R_k)
\cap (R_{\infty} \circ \theta^k) \big| \calF_k \vee \calG_\infty).
\end{equation}

Let's write $\phi (t)=-t \ln t, ~ t \geq 0$ and put
$$
I_n :=\sum_{A:A\subset S_k^{-1} R_k} \phi(\Pnum ((S_k^{-1} R_k) \cap (R_{n-k} \circ \theta^k)
=A \big| \calF_k \vee \calG_n)).
$$
Then $ H( (S_k^{-1} R_k) \cap (R_{n-k} \circ \theta^k) \big| \calF_k \vee \calG_n)= \Enum [I_n]$.

Now we estimate each item of $I_n$. Let $A\subset S_k^{-1} R_k$. Then $A$ is an $\calF_k$-measurable
random set. Since $R_n \uparrow R_\infty$,
\begin{eqnarray*}
&&\Enum\left(\left|\Pnum \Bigl( (S_k^{-1} R_k) \cap (R_{n-k} \circ \theta^k)=A \big| \calF_k
\vee \calG_n \Bigr)
-\Pnum \Bigl( (S_k^{-1} R_k) \cap (R_{\infty} \circ \theta^k)=A \big| \calF_k \vee \calG_n \Bigr)
\right| \right)\\
&\le&\Enum \left[ \Pnum \Bigl((S_k^{-1} R_k) \cap (R_{n-k} \circ \theta^k) \neq (S_k^{-1} R_k)
\cap (R_{\infty} \circ \theta^k) \big| \calF_k \vee \calG_n \Bigr) \right]\\
&=& \Pnum \Bigl((S_k^{-1} R_k) \cap (R_{n-k} \circ \theta^k) \neq (S_k^{-1} R_k) \cap (R_{\infty}
\circ \theta^k) \Bigr) \to 0
\end{eqnarray*}
Therefore,
$$
\Pnum \Bigl((S_k^{-1} R_k) \cap (R_{n-k} \circ \theta^k)=A \big| \calF_k \vee \calG_n) - \Pnum
((S_k^{-1} R_k) \cap (R_{\infty} \circ \theta^k)=A \big| \calF_k \vee \calG_n)\xrightarrow{L_1}0.
$$
Furthermore, we apply the backward martingale convergence theorem to get
$$
\Pnum ((S_k^{-1} R_k) \cap (R_{\infty} \circ \theta^k)=A \big| \calF_k \vee \calG_n) \to \Pnum
((S_k^{-1} R_k) \cap (R_{\infty} \circ \theta^k)=A \big| \calF_k \vee \calG_\infty),~a.s..
$$
Hence  for each $A\subset S_k^{-1}R_k$,
$$
\Pnum ((S_k^{-1} R_k) \cap (R_{n-k} \circ \theta^k)=A \big| \calF_k \vee \calG_n) \rightarrow
\Pnum ((S_k^{-1} R_k) \cap (R_{\infty} \circ \theta^k)=A \big| \calF_k \vee \calG_\infty).
$$
in probability.
On the other hand, since  there are  at most $2^{k+1}$ different subsets of $S_k^{-1}R_k$, we
have $I_n\le \ln(2^{k+1})=(k+1)\ln 2$ and
$$
I_n \rightarrow  \sum_{A\subset S_k^{-1}R_k} \phi(\Pnum ((S_k^{-1} R_k) \cap R_{n-k} \circ
\theta^k=A \big| \calF_k \vee \calG_\infty))
$$
in probability. Then the Dominated Convergence Theorem tells us
\begin{align*}
\lim_n\Enum [I_n] =& \Enum [\sum_{A\subset S_k^{-1}R_k} \phi(\Pnum ((S_k^{-1} R_k) \cap R_{n-k}
\circ \theta^k=A \big| \calF_k \vee \calG_\infty))]\\
=&H( (S_k^{-1} R_k) \cap R_{\infty} \circ \theta^k \big| \calF_k \vee \calG_\infty).
\end{align*}
So we have (\ref{e:1130}) and finish the proof of the lemma.
\qed

\section{Proof of Theorem \ref{thm: 1.2}}\label{sec: 4}
In order to prove Theorem \ref{thm: 1.2} we need the following result.
\begin{lem}\label{l:eneqv}
Let $\tilde S_n=X_1^{-1}X_2^{-1}\cdots X_n^{-1}$ for all $n$. Let $\tilde R_n$ be the range of $\tilde S$
until time $n$. Then
(a) $\tilde S$ is recurrent if and only if $S$ is recurrent;
(b)$\lim\limits_n \frac{H(\tilde R_n)}{n}=\lim\limits_n \frac{H(  R_n)}{n}.$
\end{lem}
\Proof
(a) Suppose that $S$ is recurrent. Then
$$
\sum_{n}\Pnum(S_i\not=e {\rm ~for~any~}1\le i< n {\rm~and~} S_n=e)=1.
$$
Since $(X_1,\cdots,X_n)$ has the same distribution with $(X_n,\cdots,X_1)$, we have
\begin{align*}
&\Pnum(S_i\not=e {\rm ~for~any~}1\le i< n, {\rm~and~} S_n=e)\\
 =&\Pnum(X_1\cdots X_i\not=e {\rm ~for~any~}1\le i<n, {\rm~and~} X_1\cdots
 X_n=e)\\
 =&\Pnum(X_n\cdots X_{n-i+1}\not=e {\rm ~for~any~}1\le i<n, {\rm~and~} X_n\cdots
 X_1=e)\\
 =&\Pnum(X_{n-i+1}^{-1}\cdots X_{n}^{-1}\not=e {\rm ~for~any~}1\le i< n, {\rm~and~} X_1^{-1}\cdots
 X_n^{-1}=e)\\
 =&\Pnum(X_{1}^{-1}\cdots X_{n-i}^{-1}\not=e {\rm ~for~any~}1\le i< n, {\rm~and~} X_1^{-1}\cdots
 X_n^{-1}=e)\\
 =&\Pnum(\tilde S_i\not=e {\rm ~for~any~}1\le i<n, {\rm~and~} \tilde
 S_n=e).
\end{align*}
So, we have $ \sum_{n} \Pnum (\tilde S_i\not=e {\rm ~for~any~}1\le i<n {\rm~and~} \tilde S_n=e)=1,$
which implies $\tilde S$ is recurrent.

(b) Since $(X_1,\cdots,X_n)$ has the same distribution with $(X_n,\cdots,X_1)$, we have
\begin{align*}
H(R_n, S_n)\ge& H(S_n^{-1}\cdot R_n)\\
=&H(\{X_n^{-1}\cdots X_i^{-1}: 1\le i\le n\}\cup\{e\})\\
=&H(\{X_1^{-1}\cdots X_{n-i+1}^{-1}: 1\le i\le n\}\cup\{e\})=H(\tilde R_n).
\end{align*}
By Theorem \ref{thm: 1.1}, we have
$$
\lim_n \frac{H(  R_n)}{n}=\lim_n \frac{H(  R_n, S_n)}{n}\ge \lim_n \frac{H( \tilde R_n)}{n}.
$$
Similarly, we have $\lim_n \frac{H( \tilde R_n)}{n}\ge \lim_n \frac{H(  R_n)}{n}$, and so we have (b).
\qed

\noindent{\bf Proof of Theorem \ref{thm: 1.2}.}
We shall prove $(1)\Rightarrow (2)\Rightarrow (3)\Rightarrow (4)\Rightarrow (1)$.\\
 {\bf Part Proof:  (1)$\Rightarrow$(2).} Let $\tilde \mu$ be the probability measure on $G$ such that $\tilde
\mu(x)=\mu(x^{-1})$ for all $x\in G$. Define $\tilde S$ and $\tilde
R_n$ as in Lemma \ref{l:eneqv}. Then  $\tilde S$ is a $(G, \tilde \mu)-$random walk.
 For each $x\in G$, set
$$
\tau_x=\inf\{n\ge 1:  S_n=x\}~~~{\rm and~~~} \tilde\tau_x
:=\inf\{n\ge 1: \tilde S_n=x\}.
$$
Write $G_{\mu}=\{x: \Pnum( \tau_{x}<\infty)=1\}$ and $G_{\tilde \mu}=\{x: \Pnum(\tilde \tau_{x}<\infty)=1\}$.
By Lemma \ref{l:eneqv}(b), to prove  (1) implies (2), it suffices to show that (1) implies $G_{\tilde\mu}
\not=\emptyset$. Since $H(X_1)>0$, we have $|{\supp}(\mu)|\ge 2$.  The rest proof of (1)$\Rightarrow$(2)
is divided into three cases.

Case I: there is an element $a\in G\backslash\{e\}$ such that ${\supp}(\mu)=\{a,e\}$. Then
${\supp}(\tilde\mu)=\{a^{-1},e\}$. It implies that $\Pnum(\tilde\tau_{a^{-1}}<\infty)=1$ and
$G_{\tilde\mu}\not=\emptyset$ as desired.

Case II: there are two different elements $a,b\in {\supp}(\mu)\cap G_{ \mu}$. Then
\begin{equation}\label{e:1312}
0=\Pnum( \tau_{a}=\infty)\ge \Pnum( S_1=b)\Pnum(
\tau_{a}=\infty| S_1=b)= \mu(b)\Pnum(
\tau_{b^{-1}a }=\infty),
\end{equation}
which implies $b^{-1}a \in G_{ \mu}$. Similarly, $a^{-1}b \in G_{\mu}$. Hence $e=(b^{-1}a)(a^{-1}b)
\in G_{\mu}$, which implies $S$ is recurrent. By Lemma \ref{l:eneqv}, process $\tilde S$ is recurrent
and so $G_{\tilde\mu}\not=\emptyset$.

Case III: there is an element $g\in G\backslash\{e\}$ such that $0<\mu(g)<1$ and $P(\tau_g=\infty)>0$.
Let us show the following inequality first
\begin{equation}\label{e:1215}
h_{ R}( \mu)\ge -\Pnum(
\tau_g=\infty)\Pnum(\tilde\tau_{e}=\tilde\tau_{g}=\infty)\ln[1- \mu(g)].
\end{equation}
To prove the statement above, for  $A\subset G$ we set
$$
\partial A:=\partial_g A=\{x\in A: x  g\not\in A\}.
$$
Similarly as \cite[Lemma 3 and Corollary 4]{BKYY2010}, one can get $\Pnum ( R_n=A)\le \left[1- \mu(g)\right]
^{|\partial A|-1}$ and
\begin{equation}\label{e:1521}
H( R_n)\ge -(\Enum|\partial  R_n|-1)\ln [1- \mu(g)].
\end{equation}
So we need estimate $|\partial  R_n|$ further.  For $x \in G\backslash\{e,g^{-1}\}$ and $ 1\le m\le n$, set
$$
F_{x,m}:=\{ S_i\not\in \{x,x  g\} {\rm~for~all~} 1\le i<m\}\cap \{ S_m=x\}\cap\{ S_j\not=x g {\rm~for~all~}
m<j \le n\}.
$$
Then $\{x\in \partial  R_n\}=F_{x,1}\cup\cdots\cup F_{x,n}$. Hence
$$
\Enum (|\partial   R_n|)\ge \sum_{1\le m\le n}\sum_{x\not=e,g^{-1}} \Pnum (F_{x,m}).
$$
Furthermore, by the Markov property
\begin{equation*}
\Pnum (F_{x,m})= \Pnum (  S_i\not\in \{x,x  g\},\forall 1\le i<m, S_m=x)\Pnum(  \tau_g>n-m).
\end{equation*}
Since $(X_1,\cdots,X_m)$ has the same law with $(X_m,\cdots, X_1)$,
we have
\beq
&&\sum_{x\not=e,g^{-1}}\Pnum (S_i\not\in \{x,x  g\},\forall 1\le i<m, S_m=x)\\
&=&\sum_{x\not=e,g^{-1}}\Pnum (X_1X_2\cdots X_i\not\in \{x,x g\},\forall 1\le i<m, X_1X_2\cdots X_m=x)\\
&=&\sum_{x\not=e,g^{-1}}\Pnum (X_m X_{m-1}\cdots X_{m-i+1}\not\in \{x,x g\},\forall 1\le i<m, X_m X_{m-1}
\cdots X_1=x)\\
&=&\sum_{x\not=e,g^{-1}}\Pnum (X_1^{-1}X_2^{-1}\cdots X_{m-i}^{-1}\not\in\{e,g\},\forall 1\le i<m,
X_1^{-1}X_2^{-1} \cdots X_m^{-1}=x^{-1})\\
&=& \Pnum (X_1^{-1}X_2^{-1}\cdots X_{i}^{-1}\not\in\{e,g\},\forall 1 \le i\le m)=\Pnum(\tilde S_i
\not\in \{e,g\}, \forall 1\le i\le m).
\eeq
Hence
\begin{eqnarray*}
\Enum (|\partial  R_n|) &\ge& \sum_{1\le m\le n}\Pnum( \tau_g>n-m) \cdot \Pnum(\tilde S_i\not\in \{e,g\},
\forall 1\le i\le m)\\
&\ge& n\Pnum ( \tau_g=\infty) \cdot \Pnum (\tilde{\tau}_e=\tilde{\tau}_g=\infty).
\end{eqnarray*}
Therefore,  (\ref{e:1215}) holds.

Now assume condition (1). Applying (\ref{e:1215}), we then get $\Pnum (\tilde{\tau}_e=\tilde{\tau}_g
=\infty)=0$, which implies $\Pnum ({\tilde\tau}_{e}<{\tilde\tau}_{g}) +\Pnum ({\tilde\tau}_{g}
<{\tilde\tau}_{e})=1$. By the strong Markov property,
\begin{eqnarray*}
\Pnum( {\tilde\tau}_{g}<\infty) &=& \Pnum( {\tilde\tau}_{g}< {\tilde\tau}_{e})+\Pnum( {\tilde\tau}_{e}
< {\tilde\tau}_{g}, {\tilde\tau}_{g}<\infty)\\
&=&\Pnum( {\tilde\tau}_{g}< {\tilde\tau}_{e})+\Pnum( {\tilde\tau}_{e}< {\tilde\tau}_{g}) \cdot \Pnum
({\tilde\tau}_{g}<\infty).
\end{eqnarray*}
So either $\Pnum( {\tilde\tau}_{e}< {\tilde\tau}_{g})=1$ or $\Pnum( {\tilde\tau}_{g}<\infty)=1$ holds,
which proves ${\bf (1)}\Rightarrow {\bf (2)}$.
\qed

\noindent{\bf Part Proof:  (2)$\Rightarrow$(3).}
Write $G_\mu=\{x:\Pnum(\tau_x<\infty)=1\}$, and set
$$
\tau_x :=\inf\{n\ge 1: S_n=x\}~~{\rm and~}~ T_x:=\inf\{n\ge 0: S_n=x\}~~{\rm for~ }x\in G.
$$
By the condition {\rm (2)}, we have $G_\mu\not=\emptyset$. If $|G_\mu\cap \supp(\mu)|\ge 2$, then as Case
II in ${\bf (1)}\Rightarrow {\bf (2)}$, we can prove $S$ is recurrent. If $e\in G_\mu$ then $S$ is recurrent,
too. So, we may assume that $S$ is transient, $|G_\mu\cap \supp(\mu)|\le 1$ and there exists $x_0\in
G_\mu$ such that $x_0\not=e$ in the following.

If $y\in G$ with $\Pnum(T_y<\tau_{x_0})>0$, then $y^{-1}x_0\in G_\mu$. This is because
$$
\Pnum(\tau_{x_0}=\infty) \ge \Pnum(T_y<\tau_{x_0}) \cdot \Pnum^y(\tau_{x_0}=\infty)
=\Pnum(T_y<\tau_{x_0}) \cdot \Pnum(\tau_{y^{-1}x_0}=\infty).
$$
For any $a\in G$ with $\Pnum(S_{\tau_{x_0}-1}=x_0a)>0$, we know
$$
\Pnum(T_{x_0a}<\tau_{x_0})>0 \hbox{ and } \Pnum(X_1=a^{-1})>0,
$$
which implies $a^{-1}\in G_\mu\cap {\supp} (\mu)$. This argument means that for any
$x\in G_\mu\backslash \{e\}$ and $b\in G$ with $\Pnum(S_{\tau_{x}-1}=xb)>0$, we have
$b^{-1}\in G_\mu\cap {\supp}(\mu)$. By the assumption $|G_\mu\cap \supp(\mu)|\le 1$, we
have $a=b$ and so,
$$
\Pnum(S_{\tau_{x}-1}=xa )=1~~{\rm for~all~}x\in G_\mu.
$$

Further, if $x\in G_\mu\setminus \{e,a^{-1},a^{-2},\cdots\}$ then $\{S_{\tau_{x}-1}=xa \}\subseteq\{\tau_{xa }
+1 \le  \tau_{x}<\infty\}$  and so, $\tau_{xa }+1\le \tau_{x}<\infty~~a.s.$. Iterating the inequality, we
get $\tau_x\ge n+\tau_{xa^{ n}} ~~a.s. $ for any $n\ge 1$, which contradict to $\Pnum(\tau_x<\infty)=1$. Hence
$$
G_\mu\subset \{e,a^{-1},a^{-2},\cdots\}.
$$
Since  $\Pnum(T_x<\tau_{a^{-1}})\ge \mu(x)>0$ for $x\in {\supp}(\mu)\setminus\{a^{-1}\}$, as before we have
$$
x^{-1}a^{-1}\in G_\mu{\rm~~~for~} x\in {\supp}(\mu)\setminus\{a^{-1}\}.
$$
Therefore,
$$
{\supp}(\mu)\subset \{a^{-1}b^{-1}: b \in G_\mu\}\cup\{a^{-1}\}\subset \{a^{-1}, e, a, a^2,\cdots\}.
$$
So, we can view $S$ as a simple random walk  on $\mathbb{Z}$ without left jump. Since $S$ is transient and
$G_\mu\subset \{e,a^{-1},a^{-2},\cdots\}$, process $S$ must be a random walk escaping to $-\infty$ without
left jump. Anyway, {\rm (3)} is true.
\qed

\noindent{\bf Part Proof: (3)$\Rightarrow$(4).}
If $S$ is recurrent, then $R_{\infty}=R_{\infty}\circ \theta=G$ almost surely. Similar as the proof of
case 1 of  Theorem 6.7.3 in \cite{Durrett2010}, one can show that $\calG_\infty$ is trivial and hence
(4) holds for this case.

Now we assume that $S$ is a random walk escaping to $-\infty$ without left jump. Then  $G$ is isomorphic
to $\mathbb{Z}$. As we are just concerning the entropies and the tail $\sigma-$algebra,  we may assume
further that $\mu$ is supported on $\{k\in \mathbb{Z}: k\ge -1\}$ with $\Enum X_1\in (-1,0)$ (and hence
$\lim_n S_n=-\infty$ a.s.).

Let $\eta:=\sup\{S_n:n\ge 0\}$ and $\xi:=\eta\circ \theta$. Then $\eta,\xi<\infty$ a.s.. Note that
$$
R_\infty=(-\infty,\eta]\cap\mathbb Z.
$$
So in order to prove the first equality of {\rm (4)}, it suffices to show $H(\xi|\eta)=H(\eta|\xi)$.
Write $\eta_n=\eta\wedge n, \xi_n=\xi\wedge n$. Since $\eta$ and $\xi$ have the same distribution,
$$
H(\eta_n|\xi_n)=H(\xi_n,\eta_n)-H(\xi_n)=H(\xi_n,\eta_n)-H(\eta_n)=H(\xi_n|\eta_n).
$$
So, it also suffices to show that  $H(\eta_n|\xi_n)\to H(\eta|\xi)$ and $H(\xi_n|\eta_n)\to H(\xi|\eta)$
as $n\to\infty$.

Put $\phi (t):=-t \ln , \forall t \geq 0$. A direct calculation gives
\beq
H(\eta_n|\xi_n)&=&\sum_{i=0}^{n-1} \Pnum(\xi=i)\sum_{j=0}^{n-1} \cdot \phi (\Pnum(\eta=j|\xi=i))\\
&&+\sum_{i=0}^{n-1} \Pnum(\xi=i) \cdot \phi (\Pnum(\eta\ge n|\xi=i))\\
&&+\Pnum(\xi\ge n) \cdot H(\eta_n|\xi\ge n)\\
&=:&I_1+I_2+I_3. \eeq By the Monotone Convergence Theorem, we have
$\lim_n I_1= H(\eta|\xi)$.  By the Dominated Convergence Theorem,
 $\lim_nI_2=0$. Note that
$\eta=\max\{0,X_1,X_1+\xi\}$. We have $\eta\ge \xi-1$ and hence
$\eta_n\in \{n-1,n\}$ if $\xi\ge n$. So, $$H(\eta_n|\xi\ge n)\le \ln
2.$$ It follows  $\lim_nI_3= 0$ immediately. Therefore
$\lim_{n\to\infty}H(\eta_n|\xi_n)=H(\eta|\xi)$.

Similarly, in order to prove $\lim_{n\to\infty}H(\xi_n|\eta_n)=H(\xi|\eta)$, we need only to show that
\begin{equation}\label{e:1546}I_3':=\Pnum(\eta\ge n)H(\xi_n|\eta\ge n)\to
0.\end{equation}Since   $\xi$ and $X_1$ are independent, and
$\eta=\max\{0,X_1,X_1+\xi\}$, we have
\begin{align*}
&\sum_{j=0}^{n-1} \phi (\Pnum(\xi=j, \eta \ge n)) = \sum_{j=0}^{n-1} \phi (\Pnum(\xi=j, X_1\ge n-j))\\
=&\sum_{j=0}^{n-1} \phi (\Pnum(\xi=j)) \cdot \Pnum( X_1\ge n-j) +\sum_{k=1}^{n} \Pnum( X_1\ge k) \cdot \phi
(\Pnum (\xi=n-k)).
\end{align*}
By the Dominated Convergence Theorem and the $L^1$-integrability of $X_1$, we get that both the two summations
above tend to zero as $n\rightarrow\infty$. Hence
\begin{align*}
I_3'\le \phi (\Pnum(\xi\ge n, \eta\ge n)) +\sum_{j=0}^{n-1} \phi (\Pnum(\xi=j, \eta\ge n))\rightarrow 0.
\end{align*}
As $I_3'\ge 0$ for all $n$, we then have (\ref{e:1546}). This proves the first equality of {\rm (4)}.

We proceed to prove $\mathcal{G}_\infty\equiv\sigma(R_\infty)\mod\Pnum$.
Clearly, $\calG_\infty\supseteq\sigma(R_\infty)$. Since $R_\infty=(-\infty,\eta]\cap \mathbb Z$, we need only
to show $\calG_\infty\subseteq\sigma(\eta)\mod\Pnum$. Write
$$
\rho=\inf\big\{n\ge 0: S_n=0,~\sup\{S_k:k\ge n\}\le \max\{S_k:k\le n\}~\big\}.
$$
As $\lim_n S_n=-\infty$ a.s., we have $\rho<\infty$ a.s.. For $n\ge 0$, set $M_n=\inf\{S_k:k\le n\}$. Then
$R_n=[ M_n,\eta ]\cap\mathbb Z$ conditioned on $\rho\le n$. For $E\in\calG_\infty$ and $n\ge 1$, there exists
non-random set $A_n$ such that
$$
E=\big\{\{(R_k,S_k): k\ge n\} \in A_n\big\}.
$$
Let $E_n=\big\{\{([M_k,\eta], S_k): k\ge n\} \in A_n\big\}$. Then conditioned on $\{\rho<\infty\}$, there has
$$
E= \bigcup\limits_{n\ge 1}^\infty\bigcap\limits_{k=n}^\infty E_n\in \mathcal{H},
$$
where $\mathcal{H}:=\bigcap_n\sigma(\eta, M_k,S_k:k\ge n)$. So, from  $\rho<\infty$ a.s. we get
$$
\calG_\infty\subseteq\mathcal{H} \mod\Pnum.
$$
As the possible values of $\eta$ are countable,
$$
\mathcal{H}=\sigma(\eta)\vee\bigcap_n\sigma(M_k,S_k:k\ge n).
$$
Let $T_0=0$. Set $T_{k}=\inf\{n> T_{k-1}: X_n=-k\}$ and $U_k=(S_{T_{k-1}},\cdots, S_{T_{k}})+k$ for $k\ge 1$.
Then $U_1,U_2,\cdots$ are i.i.d., and $\bigcap_n\sigma(M_k,S_k:k\ge n)$ is contained in the exchangeable
field of  $U_1,U_2,\cdots$. By Hewitt Savage 0-1 law, $\bigcap_n\sigma(M_k,S_k:k\ge n)$ is trivial. Therefore
we have $\calG_\infty\subseteq \mathcal{H} \equiv \sigma(\eta)\mod\Pnum$ as desired.
\qed

{\bf Part Proof: (4)$\Rightarrow$(1).}
By {\rm (4)}, we get
\begin{eqnarray*}
H(X_1,R_{\infty}\circ\theta|\calG_\infty)&=&
H(X_1,R_{\infty}\circ\theta|R_\infty)=H(R_{\infty}\circ\theta|R_\infty)+H(X_1|R_{\infty}\circ\theta,R_\infty)
\\
&=&H(R_{\infty}|R_\infty\circ\theta)+H(X_1|R_{\infty}\circ\theta,R_\infty)\\
&=&H(X_1,R_\infty|R_\infty\circ\theta)=H(X_1|R_\infty\circ\theta)=H(X_1).
\end{eqnarray*} By Theorem \ref{thm: 1.4'}, it follows immediately
$$
h_R(\mu)=H(X_1)-H (X_1, R_\infty\circ\theta \big| \calG_\infty)=0.
$$
We have completed the proof of Theorem \ref{thm: 1.2}.
\qed

\section{Proof of Theorem \ref{thm: 1.3}}\label{sec: 5}

In order to prove this theorem, we need more notations. We encode each realization of the random weighted
graph $\Gamma_n$ as the following. Set $G=\{g_i: i\ge 0\}$ with $g_0=e$ and write $\mathbb{V}
=\cup_{l=1}^\infty\mathbb{N}^l$. We give a dictionary order of $\mathbb{V}$. That is, for any two vector
$x\in \mathbb{N}^i$ and $y=\mathbb{N}^j$, we say that $x$  prior to $y$
if either $i<j$ or $i=j$ with $x^{(1)}=y^{(1)},\cdots, x^{(l-1)}=y^{(l-1)},
x^{(l)}<y^{(l)}$ for some $l\le i$. Define
$\zeta:\mathbb{V}\rightarrow G$ as follows:
$$
\zeta(v)=g_{v^{(1)}}\cdots g_{v^{(l)}}~~{\rm~for~}v=(v^{(1)},\cdots,
v^{(l)})\in \mathbb{N}^l {\rm ~and~}l\ge 1.
$$

Now fix $(A,B,C)$ an arbitrary realization of $\Gamma_n$. We   use the width
first research algorithm to visit every vertex of the digraph $(A,B)$.

Step 0: visit  $v_0=e$. Set $N_0=\{i\ge 0:(e,g_i)\in B\}$ and $V_0=N_0\setminus\{0\}$.

$\cdots$

Step $h< |A|$:  choose $x\in V_{h-1}$, so that $x$ prior to all other elements of $V_{h-1}$;
visit $v_h=\zeta(x)$; and set $N_h=\{i\ge 0: (v_h,v_hg_i)\in  B\}$ and
$$
V_h=(V_{h-1}\cup\{(x,i): i\in  N_h\})\setminus \{v\in \mathbb{V}:\zeta(v)\in \{v_0,\cdots, v_h\}\}.
$$

$\cdots$\\
Since $\Gamma_n$ is generated by the first $n$ steps of random walk $S$ and $\mathbb{P}( \Gamma_n
=(A,B,C))>0$,  for each $v\in A\setminus\{e\}$ there exists at least one oriented path from $e$ to
$v$. This ensures that at each step of the  algorithm above, we always have $V_{h-1}\not=\emptyset$.
Furthermore,
$$
\{v_h: 0\le h<|A|\}=A{\rm~and~}  ~V_{|A|-1}=\emptyset.
$$
Note that  $v_0$ and $V_0$ are determined by  $N_0$. So, $v_h$ and $V_h$ are determined by $N_0,
\cdots, N_h$ for every $h< |A|$. Especially, $v_0,\cdots,v_{|A|-1}$ are determined by $N_0,
\cdots, N_{|A|-1}$. Therefore, $\{N_h:  h<|A|\}$ are one-to-one corresponding to $(A,B)$.

Next, for each $i,h\ge 0$, set
$$
O^{i,h}=\begin{cases}
C((v_h,v_hg_{i})),&{\text~if~}h<|A| \text {~and~} (v_h,v_hg_i)\in B;\\
0,&{\text~otherwise}.
\end{cases}
$$
Then for $h< |A|$ we have $ N_h=\{i\ge 0: O^{i,h}>0\}$, which implies $(A,B,C) \mapsto \{O^{i,h}:
i,h\ge 0\}$  is also one to one.

Now we  define a sequence of random variables. For  $i,h\ge 0$ and $n\ge 1$, set
$$
O^{i,h}_n=O^{i,h}(\Gamma_n).
$$
Write $O_n=\{O^{i,h}_n: i,h\ge 0\}$. As the statement above, $\Gamma_n\mapsto O_n$ is one to one,
and so
\begin{equation}\label{e:1935}
H(\Gamma_n)= H(O_n).
\end{equation}
Set $O_n^i=\sum_{h=0}^\infty O_n^{i,h}$. Then $O_n^i=\sum_{k=1}^n 1_{\{X_k=g_i\}}$, which implies
$O_n^i\sim B(n,\mu(g_i))$ and
\begin{equation}\label{e:1936}
\lim_n\frac{O_n^i}{n}=\mu(g_i)~~{\rm~ a.s.}.
\end{equation}
The following result plays a key role in the upper bound of $H(\Gamma_n)$.

\begin{lem}\label{l:2127}
Let $Y_n^i=O_n^i\ln(1+\frac{|R_n|-1}{O_n^i})+(|R_n|-1)\ln(1+\frac{O_n^i}{|R_n|-1}).$
Then
$$
\varlimsup_n \frac{H(\Gamma_n)}{n}\le \varlimsup_n \frac{1}{n} \sum_{i=0}^\infty \Enum \left(Y_n^i\right).
$$
\end{lem}
\Proof
By (\ref{e:1935}),
\begin{align*}
H(\Gamma_n)= H(O_n) \le& H(|R_n|)+H(O_n^i, i\ge 0)+H(O_n\big|~|R_n|, O_n^i, i\ge 0).
\end{align*}
We estimate the three terms separately. First,  $\lim_n n^{-1}H(|R_n|)=0$ since $|R_n|\le n+1$. Next, fix
$M\in \Nnum$ and calculate
\begin{align*}
H(O_n^i,i\ge 0)\le& \sum_{i=0}^{M-1} H(O_n^{i})+\sum_{i=1}^n H(X_i1_{\{X_i\not\in \{g_j: j< M\}\}})\\
\le& M\ln (n+1)+n H(X_11_{\{X_1\not\in \{g_j: j< M\}\}}).
\end{align*}
Hence
$$
\lim_n n^{-1}H(O_n^i,i\ge 0)\le H(X_11_{\{X_1\in\{g_j:j\ge M\}\}}).
$$
By the arbitrariness of $M$ and noting $H(X_1)<\infty$,  we get that the right side of the above inequality
vanishes as $n$ goes to infinity.

We shall estimate the third term. Let $k\in \mathbb{N}$ and $a_i\in \mathbb{Z}^+,i\ge 0$ be constants. If
$|R_n|=k$ and  $O_n^i=a_i$ for  $i\ge 0$, then for each $i\ge 0$, there are $\binom{a_i+k-1}{k-1}$
possibilities at most for the choices of $(O_n^{i,h},h\ge 0)$ since $\sum_{h=0}^{k-1}O_n^{i,h}=a_i$. Hence
the number of different choices of $O_n$ is $\prod_{i\ge 0} \binom{a_i+k-1}{k-1}$ at most. Therefore,
$$
H(O_n\big| |R_n|=k, O_n^i=a_i,i\ge 0)\le \sum_{i=0}^\infty\ln \binom {a_i+k-1}{k-1}.
$$
From the Stirling formula, one knows
$$
\binom{m}{l}\le \left(\frac{m}{l}\right)^l \left(\frac{m}{m-l}\right)^{m-l}~{\rm~for~all~}m>l\ge 1.
$$
Therefore, for $k\in \mathbb{N}$ and $\{a_i: i\ge 0\}\subset \mathbb{Z}^+$ we have
\begin{align*}
H(O_n\big| |R_n|=k, O_n^i=a_i,i\ge 0) \le& \sum_{i=0}^\infty\left(a_i\ln(1+\frac{k-1}{a_i}) +(k-1) \ln
(1+\frac{a_i}{k-1})\right),
\end{align*}
which implies
\begin{align*}
H(O_n\big|~ |R_n|, O_n^i, i\ge 0)\le \sum_{i=0}^\infty\Enum (Y_n^i).
\end{align*}
Combining the estimates of the three  terms, we get the desired result.
\qed

\noindent\textbf{Proof of  Theorem \ref{thm: 1.3}.}
Suppose that $S$ is transient. Let $\tau_x=\inf\{n\ge 1:S_n=x\}$ and $T_x=\inf\{n\ge 0:S_n=x\}$ for all
$x\in G$, where $\inf\emptyset=\infty$. Then $\Pnum(\tau_e<\infty)<1$. Since $\Pnum(\tau_e=\infty)
=\sum_{a\not= e}\mu(a)\Pnum(\tau_{a^{-1}}=\infty)$, there exists some $a\in G\backslash \{e\}$ such that
$\mu(a)\in (0,1)$ and $\Pnum(\tau_{a^{-1}}=\infty)>0$. For each realization  $(A,B,C)$  of $\Gamma_n$, set
$$
\diamondsuit_a (A,B,C)=\big\{x\in A: \{g\in G:(x,xg)\in B\}=\{a\}\big\}.
$$
By the strong Markov property $$\Pnum(\Gamma_n=(A,B,C))\le \Pnum(X_{T_x+1}=a {\rm~for~all~} x \in
\diamondsuit_a (A,B,C)   )\le \mu(a)^{|\diamondsuit_a(A,B,C)|},
$$
So,
$$
H(\Gamma_n) \ge -\Enum|\diamondsuit_a\Gamma_n|~\ln [\mu(a)].
$$
Note that if $T_x<n$, $X_{T_x+1}=a$ and $S_{k}\not=x$ for all $k>T_x+1$, then $x \in \diamondsuit_a\Gamma_n$.
Hence
$$
\Enum| \diamondsuit_a\Gamma_n|\ge \sum_x \Pnum(T_x<n)\mu(a) \Pnum(\tau_{a^{-1}}=\infty) =\mu(a) \Pnum
(\tau_{a^{-1}}=\infty)\Enum|R_{n-1}|.
$$
Since $S_n$ is transient, $\lim \limits_{n\to\infty} \Enum|R_n|/n=\gamma_{\rm escape} :=\Pnum (S_n \neq e,
\forall n \geq 1)$ by \cite{Derriennic1983}. Therefore
$$
\lim_{n\to\infty} \frac{H(\Gamma_n)}{n} \ge [-\mu(a) \cdot \ln \mu(a)] \cdot \Pnum(\tau_{a^{-1}}=\infty)
\cdot \gamma_{\rm escape}>0.
$$

On the other hand, suppose $S$ is recurrent. Then by \cite{Derriennic1983} we have $|R_n|/n\rightarrow 0~a.s.$.
Recall the definition of $Y_n^i$ in Lemma \ref{l:2127}. Using (\ref{e:1936}), we conclude that with probability
one
$$
\lim_n \frac{Y_n^i}{n}\le \lim_n\frac{|R_n|-1}{n} \ln (1+\frac{n}{|R_n|-1})+\frac{O_n^i}{n}\ln
(1+\frac{|R_n|-1}{O_n^i})=0.
$$
Note that $Y_n^i\le  (|R_n|-1)+O_n^i\le 2 n$ since $\ln(1+x) \le x$ for all $x\ge 0$. Then by the Dominated
Convergence Theorem,
\begin{equation}\label{e:2111}
\lim_n\frac{\Enum (Y_n^{i})}{n}=0.
\end{equation}
On the other hand, using  $\ln(1+x) \le x$ again we have
$$
Y_n^i\le O_n^i\ln\left(\frac{2n}{O_n^i}\right)+O_n^i\le O_n^i (2+\ln n)-O_n^i\ln O_n^i.
$$
Since $f(x)=-x\ln x$ is concave on $(0,\infty)$ and $O_n^i\sim B(n, \mu(g_i))$,
$$
\Enum (-O_n^i\ln O_n^i)\le -\Enum (O_n^i)\ln(\Enum (O_n^i))=-n\mu(g_i)\ln (n\mu(g_i)).
$$
So,
\begin{equation}\label{e:co17}
n^{-1}\Enum (Y_n^i)\le \mu(g_i)(2+\ln n)-\mu(g_i) \ln (n\mu(g_i))=2\mu(g_i)-\mu(g_i)\ln(\mu(g_i)).
\end{equation}
Note that
$$
\sum_{i=0}^\infty \big[2\mu(g_i)-\mu(g_i)\ln(\mu(g_i))\big]=2+H(X_1)<\infty.
$$
By  (\ref{e:2111}) and the Dominated Convergence Theorem, we obtain
$$
\lim_n\sum_{i=0}^\infty \frac{\Enum (Y_n^{i})}{n}=0.
$$
Applying  Lemma \ref{l:2127}, we get $\lim_n\frac{H(\Gamma_n)}{n}=0$  and complete the proof   of the theorem.
\qed
\section{Some Related Results}\label{sec: 6}

\subsection{Extreme Case  $h_R (\mu)=H (X_1)$}
As is known, $0 \leq h_R (\mu) \leq H (X_1)$. In this subsection,  we would  consider the  extreme case $h_R (\mu)
=H (X_1)$.
\begin{thm}\label{thm: 1.6}
$h_R (\mu)=H (X_1)$ if and only if the escape rate of process $S$ satisfies
$$
\gamma_{\rm escape}:=\Pnum (S_n \neq e, \forall n \geq 1) =1.
$$
In addition, if $h_R (\mu)=H (X_1)$ then $\sigma (R_n)=\sigma (X_k: 1 \leq k \leq n)$ and $H (R_n)=H
(R_n, S_n)=n H (\mu)$ for all $n \geq 1$.
\end{thm}
\noindent \textbf{Proof.}
Let $\gamma_{\rm escape}<1$. Then there exists $n\ge 1$ and
elements $a_1,\cdots,a_n\in \supp(\mu)$ such that
$a_1a_2\cdots a_n=e$. Since $H(X_1)>0$, we can find some $b\in
{\supp(\mu)}\backslash\{a_1\}$. Let
$$
\textbf{x}=( a_1 , \cdots, a_n, a_1,\cdots, a_n, b,
a_1,\cdots, a_n)~{\rm~and~}~
\textbf{y}=( a_1 , \cdots, a_n,b, a_1,\cdots, a_n,
a_1,\cdots, a_n).
$$
Then both $\bf{x}$ and $\bf{y}$  are possible realizations of $(X_1, \cdots,X_{3n+1})$. Moreover, if
$(X_1,\cdots,X_{3n+1})\in\{\bf{x},\bf{y}\}$ then $S_{3n+1}=b$ and $R_{3n+1}=\{e,a_1,a_1a_2,\cdots,a_1 a_2
\cdots a_{n-1}, b, ba_1,ba_1a_2, \cdots,ba_1\cdots a_{n-1}\}$. Hence $ H(X_1, \cdots, X_{3n+1}
\big| R_{3n+1}, S_{3n+1})>0, $ which implies
$$
H(R_{3n+1}, S_{3n+1})<H(X_1,\cdots,X_{3n+1})=(3n+1)H(X_1).
$$
Using (\ref{e:1649}), we then have
$$
h_R(\mu)\le \tfrac{H(R_{3n+1},S_{3n+1})}{3n+1}<H(X_1).
$$

Next, suppose that $\gamma_{\rm escape}=1$, this means $\Pnum (\tau_e<\infty)=0$. For any two $x,y\in G$,
we call $x\prec y$ if and only if $\Pnum(\tau_{x^{-1}y}<\infty)>0$. Then  the binary relation '$\prec$'
is a strict partial order (i.e., $x\prec y$ implies $x \neq y$ and not $y \prec x$).  By definition,
with probability one
$$
e=S_0 \prec S_1 \prec S_2 \prec \cdots.
$$
Therefore, for each $n+1\ge k\ge1$, element $S_{k-1}$ is almost surely the $k$-th least element among
$R_n$. It follows that $\sigma(R_n)=\sigma(S_{k}: 0\le k\le n)=\sigma(X_1,\cdots,X_n)$ and hence
$h_R (\mu) =H (\mu)$.
\qed
\begin{rem}
When $\gamma_{\rm escape}=1$, it is possible that $h_S (\mu)=0<h_R (\mu)$, e.g., $S$ is the
$(\mathbb{Z}^2,\mu)$-random walk with $\mu(\{(1,0)\})=\mu(\{(0,1)\})=\tfrac12$.
\end{rem}
\begin{rem}
If $\gamma_{\rm escape}=1$, then by Theorem \ref{thm: 1.6}, we have $H(R_\infty\circ\theta|R_\infty)=0$
and $H(R_\infty|R_\infty\circ\theta)=H(X_1)$. This gives an example that $H(R_\infty\circ\theta
|R_\infty)\not=H(R_\infty|R_\infty\circ\theta)$;so the first condition in  {\rm(4)} of Theorem
\ref{thm: 1.2} is nontrivial. It is easy to see that if $H (R_\infty)<\infty$ the condition
$H({R}_\infty\circ \theta|{R}_\infty)= H({R}_\infty|{R}_\infty\circ\theta)$ holds  automatically.
\end{rem}

\subsection{Extreme Case $h_R(\mu)=0$ for Transient SRW}
Let us fix  $S$ as  a ($\mathbb{Z},\mu)-$random walk escaping to $-\infty$ without left jump.
Recall $\eta=\sup\{S_n:n\ge 0\}$  and  $R_\infty=(-\infty,\eta]\cap \mathbb Z$ a.s. (see the Part
Proof: (3)$\Rightarrow$(4) in Section \ref{sec: 4}). We will discuss whether
$H(R_\infty)$ is finite or not.

First, we show a relation between the finiteness of the entropy of an  $\Nnum$-valued random variable $Y$
and the integrability of $\ln Y$.

\begin{lem}\label{lem: entropy-vs-integral}
Let $Y$ be an $\Nnum$-valued random variable. If $\Enum (\ln Y)<\infty$, then
$H (Y)<\infty$. Conversely, if the probability sequence
$\{p_n:=\Pnum (Y=n): n \geq 1\}$ is eventually decreasing in $n$,
then $H (Y)<\infty$ also implies $\Enum (\ln Y) <\infty$.
\end{lem}
\noindent{\bf Proof.}
Since $ -x \ln x$ is increasing on $(0, 1/4]$,
\begin{eqnarray*}
H (Y) &=& \sum_{n\ge 1} -p_n \ln p_n\\
&\leq& -p_1\ln p_1 -\sum_{n\ge 2:p_n< n^{-2}} p_n\ln p_n +\sum_{n\ge 2:p_n\ge n^{-2}} p_n \ln
\tfrac{1}{p_n}\\
&\leq& e^{-1}-\sum_{n\ge 2:p_n< n^{-2}} n^{-2}\ln (n^{-2})+ \sum_{n\ge 2:p_n\ge n^{-2}} p_n\ln n^2\\
&\le& C+2 \Enum (\ln Y),
\end{eqnarray*}
where $C=e^{-1}+2\sum\limits_{n\ge 2} \frac{\ln n}{n^2}<\infty$. This proves the first part of the lemma.

Conversely, assume for simplicity that $p_n$ is decreasing in $n$. Then we must have $p_n \leq {1}/{n}$,
and so
$$
H (Y)=\sum_n -p_n \ln p_n \geq \sum_n p_n \ln n=\Enum (\ln Y).
$$
This proves the second part of the lemma.
\qed
\begin{thm}\label{rem: 1.2}
Let $S$ be a ($\mathbb{Z},\mu)-$random walk escaping to $-\infty$
without left jump. Then,  $H (R_\infty)<\infty$ if and only if
$\Enum [|X_1| \ln |X_1|]<\infty$.
\end{thm}
\Proof
For each $k\ge -1$, set
$$
p_k:=\Pnum (X_1=k),~~p_{k+}:=\sum_{\ell=k}^\infty p_l,~~ p_{k++}:=\sum_{\ell=k}^\infty p_{l+}.
$$
Write $f_k:=\Pnum (\eta=k)$ for $k\ge 0$, and define $f_{k+}$ in a similar way. Put
$$
F (t) :=\sum_{n=0}^\infty f_n \cdot t^n, \quad P (t) :=\sum_{n=0}^\infty p_n \cdot t^n, \quad t \in [0, 1].
$$
Since $\eta=\max\{0, X_1, X_1+\eta \circ \theta\}$,
\begin{eqnarray*}
F (t)&=&\Enum [ \Enum (t^\eta \big| X_1)]\\
&=& \Pnum (X_1=-1) \cdot \Enum [t^{\max\{0, \eta -1\}}] +\sum_{n=0}^\infty \Pnum (X_1=n) \cdot \Enum
[t^{\eta+n}]\\
&=&q [\Pnum (\eta=0)+\sum_{n=1}^\infty \Pnum (\eta=n) \cdot t^{n-1}] +P (t) \cdot F (t)\\
&=&q[f_0+\frac{F (t) -f_0}{t}] +P (t) \cdot F (t),
\end{eqnarray*}
where we write $q=p_{-1}$ for convenience. Hence
\begin{equation}\label{e:0925}
F (t)=\frac{q \cdot f_0 (1-t)}{q+t \cdot P (t) -t}, \quad t \in (0, 1).
\end{equation}
Put
$$
\bar{P} (t):=\sum_{n=1}^\infty p_{n+} \cdot t^n.
$$
Then for $t\in(0,1)$,
$$
\bar{P} (t)=\sum_{l=1}^\infty p_l\frac{t-t^{l+1}}{1-t}=\tfrac{t}{1-t}\sum_{l=0}^\infty
p_l-\tfrac{t}{1-t}\sum_{l=0}^\infty p_l t^l=\frac{t(1-q)}{1-t}-\frac{tP(t)}{1-t}.
$$
So, (\ref{e:0925}) can be rewritten as
\begin{equation}\label{e:1338}
F (t)=\frac{q \cdot f_0}{q-\bar{P} (t)}, \quad t \in (0, 1),
\end{equation}
which implies
$$
F(t)-f_0=F(t)\bar{P} (t)/q.
$$
Comparing the coefficients of $t^n$ in the series expansions of both sides, we obtain
$$
f_n=\sum_{k=1}^n {p_{k+}} \cdot f_{n-k}/q~~{\rm for ~all~}n\ge 1.
$$
Therefore, for each $n\ge 1$,
\begin{equation}\label{e:1252}
f_n\ge p_{n+}f_0/q.
\end{equation}

Now, assume $H(R_\infty)<\infty$. Then $H(\eta)<\infty$ and so $\sum\limits_{n} -f_n\ln f_n<\infty$.
Since $S$ is a random walk escaping to $-\infty$, we have $f_0=\Pnum(\eta=0)=\Pnum(S_n\le 0
{\rm ~for ~all~}n)>0$. Applying (\ref{e:1252}) gives
$$
\sum\limits_{n} -p_{n+}\cdot\ln p_{n+}<\infty.
$$
Since $p_{n+}$ is decreasing in $n$, we apply Lemma \ref{lem: entropy-vs-integral} to get
$$
\sum_{n}^\infty p_{n+} \cdot \ln n<\infty.
$$
It follows  $\Enum [|X_1| \ln |X_1|]<\infty$ immediately.

Conversely, suppose $\Enum [|X_1| \ln |X_1|]<\infty$.  Then
$$
\sum_{k=1}^\infty\frac{ p_{k++}}{k} <\infty.
$$
By (\ref{e:1338}), we have $q-\bar{P}(1)=\frac{qf_0}{F(1)}=qf_0$. So for $t\in [0,1),$

$$
F(1)-F(t)=1-\frac{q  f_0}{q-\bar{P} (t)} =\frac {\bar{P}(1)-\bar{P}(t)}{q-\bar{P}(t)} \le
\frac{\bar{P}(1)-\bar{P}(t)}{q-\bar{P}(1)}=\frac {1}{qf_0}[{\bar{P}(1)-\bar{P}(t)}].
$$
Dividing by $1-t$ into both sides, we get
$$
\sum_{n=1}^\infty f_n\sum_{k=1}^{n}t^{k-1}=\frac{1}{qf_0}\sum_{n=1}^\infty
p_{n+1}\sum_{k=1}^{n}t^{k-1}.
$$
Hence
$$
\sum_{k=1}^\infty f_{k+}\cdot t^{k-1}\le \frac 1{qf_0}\sum_{k=1}^\infty p_{k++}\cdot
t^{k-1}.$$
Integrating of the both sides with respect to  $t$  on  $[0,1]$, we get
$$
\sum_{k=1}^\infty \frac{f_{k+}}{k}\le \frac {1}{qf_0} \sum_{k=1}^\infty \frac {p_{k++}}{k}<\infty,
$$
which implies $\Enum \ln (\eta +1)<\infty$. By Lemma \ref{lem: entropy-vs-integral}, $H(R_\infty)
=H (\eta+1)<\infty$.
\qed

\noindent{\sl \textbf{Acknowledgements} \quad} {We would like to express our gratitude for the financial
support from Shanghai Center of Mathematics during Sept. 2014-Jan. 2015 which makes the current co-work
possible. This work is partially supported by NSFC (No. 11531001, No. 10701026, No. 11271077 and No. 11371317)
and the Laboratory of Mathematics for Nonlinear Science, Fudan University.}





\bibliographystyle{elsarticle-num}



\end{document}